\magnification 1200
\overfullrule=0mm

\font\Bb=msbm10
\hsize 16true cm
\def\fin{\hfill \vrule height5pt width3pt depth2pt\vskip 6
pt}

\def\C{\hbox{\Bb C}}
\def\R{\hbox{\Bb R}}

\def\N{\hbox{\Bb N}}

\centerline {\bf ESTIMATES FOR JACOBI-SOBOLEV TYPE
ORTHOGONAL POLYNOMIALS}
\bigskip
\centerline {by}
\bigskip

\centerline { M. Alfaro\footnote*{This research was partially
supported by DGICYT PB93-0228-C02-02}, F. Marcell\'an\footnote{**}
{This research was partially supported by
DGICYT PB93-0228-C02-01}, M.L. Rezola*}

\vskip 2true cm
\bigskip
   Abstract. Let the Sobolev-type inner product $\displaystyle
\langle f,g
\rangle =\int_{\R}fgd\mu_{0}+\int_{\R} f'g'd\mu_1$  with
$\mu_0=w
+M\delta_c$, $\mu_1=N\delta_c$ where $w$ is the Jacobi
weight,  $c$ is
either 1 or $-1$ and $M,N \geq 0$. We obtain estimates and
asymptotic
properties on $[-1,1]$ for the  polynomials orthonormal with
respect to
$\langle .,. \rangle$ and their kernels. We also compare these
polynomials
with Jacobi orthonormal polynomials.
\vskip 3truecm  {\noindent \it AMS Subject Classification
(1991)}: 33C45,
42C05
 \vfill \eject 
{\noindent \bf
1. Introduction} \bigskip Recently the study of polynomials
orthogonal with
respect to a nonstandard inner product

$$\langle f,g \rangle =\int_{\R}fgd\mu_{0}+\sum_{k=1}^{m} 
\int_{\R}
f^{(k)}g^{(k)}d\mu_k$$ has attracted the interest of many
researchers. In
particular when $m=1$ and $\mu_1$ is an atomic measure
supported at a point
$c \in \R$, results concerning algebraic properties of such
polynomials and
the location of their zeros have been done (see for instance
[1]). From an
analytic point of view, the relative asymptotic behaviour of
such
polynomials when $\mu_0$ belongs to the class $M(0,1)$ has
been
accomplished in several papers ([2], [6] and [7]). This
behaviour is
considered in compact sets of $\C \setminus {\rm supp \enspace
\mu_0}$.

However, the behaviour of polynomials in ${\rm supp \enspace
\mu_0}$
remains an open question. The aim of this paper is to cover
this lack in
the literature. In fact, a first approach was given by
Marcell\'an and
Osilenker [8] when $m=1$, $d\mu_0=\chi_{[-1,1]} dx+M(
\delta_1+\delta_{-1}
)$ and $d\mu_1= N(\delta_1+\delta_{-1})$ using some previous
work by
Bavinck and Meijer ([3], [4]),  ($\delta_c$ denotes a Dirac
measure
supported at the point $c$).

In our paper, we will consider $m=1$ $$\eqalign{ d\mu_0(x)
&=(1-x)^{\alpha}(1+x)^{\beta}dx+M\delta_1(x) \cr d\mu_1(x)
&=N\delta_1(x)
\cr}$$ with $\alpha >-1$ and $\beta >-1$.

In Section 2 we present the basic tools concerning the
polynomials
orthogonal with respect to the above inner product with
special emphasis in
the case of the so-called Jacobi-Sobolev type polynomials and
some results
about Jacobi polynomials which we will need throughout the
paper.
Section 3 deals with pointwise analysis and upper bounds for
Jacobi-Sobolev
type polynomials as well as an upper bound of their uniform
norm using the
corresponding estimates for standard Jacobi polynomials.
Previously, we
study the behaviour of the coefficients which appear in their
representation in terms of Jacobi polynomials and, as a
consequence, an
estimate for them at the ends of the interval as well as an
estimate for
its first derivative are given.

Finally, in Section 4 we obtain some bounds and estimates for
the kernels
associated with the  polynomials considered above. In
particular, the
analogue of a very well known result by M\'at\'e-Nevai-Totik
concerning
Christoffel functions is deduced.

In such a way we can give a complete answer in order to
estimate the
behaviour on $[-1,1]$ of such polynomials. Notice that some of
the above
results, when $d\mu_0=wdx+M\delta_c$ where $w$ is a
generalized Jacobi
weight and $\mu_k=0$ $(k=1,...,m)$, have been obtained in
[5].

\vfill \eject

{\noindent \bf 2. Representation formulas and basic results}
\bigskip

Let $\mu$ be a positive Borel measure on ${\R}$ whose moments
are finite
and whose support is an infinite set.

We consider the inner product  $$\langle f,g \rangle
=\int_{\R}fgd\mu+M
f(c)g(c)+Nf'(c)g'(c) \qquad  M,N  \geq 0 \qquad  c \in \R
\eqno (1) $$ .  Let $p_n$ and $q_n$ be the polynomials orthonormal with
respect to the
measure $\mu$ and the inner product (1), respectively.
Denote $q_n(x)=\gamma_n x^n+...$ and $p_n(x)=k_n
x^n+...$ . The Fourier
expansion of  $q_n$ in terms of $p_k$ $(k=0,...,n)$ leads to
$$q_n(x)={\gamma_n \over
k_n}p_n(x)-Mq_n(c)K_{n-1}(x,c)-Nq_n'(c)
K_{n-1}^{(0,1)}(x,c) \eqno (2)$$
 We have used the abbreviation 
$$K_n^{(r,s)}(x,y)=\sum_{k=0}^n
p_k^{(r)}(x)p_k^{(s)}(y)={\partial^{r+s} \over \partial x^r
\partial
y^s}K_n(x,y)$$ where, as usual, $K_n(x,y)=\sum_{k=0}^n
p_k(x)p_k(y)$.

If we take derivatives in (2) with respect to $x$ and
evaluating at $x=c$,
the values of $q_n(c)$ and $q_n'(c)$ can be expressed by
$$\eqalign {
q_n(c) &={\gamma_n \over k_n D_n}
[p_n(c)\{1+NK_{n-1}^{(1,1)}(c,c)\}-Np_n'(c)K_{n-1}^{(0,1)}(c,c)]\cr
q_n'(c) &={\gamma_n \over
k_nD_n}[-Mp_n(c)K_{n-1}^{(0,1)}(c,c)+p_n'(c)\{1+MK_{n-1}(c,c)\}]\cr } $$
where $$D_n=1+MK_{n-1}(c,c)+NK_{n-1}^{(1,1)}(c,c)+MN[K_{n-1}(c,c)K
_{n-1}^{(1,1)}(c,c)
-(K_{n-1}^{(0,1)}(c,c))^2] \eqno (3)$$ ( note that
$D_n=D_n(M,N,c)>0$ for
all $M \geq 0$, $N \geq 0$ and $c \in \R$). \medskip Let
$p_n(x;\mu_j)=k_n(\mu_j)x^n+...$, $j=0,1,2,...$, the
orthonormal
polynomials with respect to the measure $d\mu_j=(x-c)^{2j}
d\mu$ (where
$\mu_0=\mu$) and $K_n(x,y;\mu_j)$ the corresponding kernels.
Expanding
$(x-c)p_{n-1}(x;\mu_{j+1})$ in terms of $p_k(x;\mu_j)$ we
obtain (see [1,
Lemma 2.1]) $$(x-c)p_{n-1}(x;\mu_{j+1})={k_{n-1}(\mu_{j+1})
\over
k_n(\mu_j)}[p_n(x;\mu_j)-{p_n(c;\mu_j) \over
K_{n-1}(c,c;\mu_j)}K_{n-1}(x,c;\mu_j)]$$ Using the
orthonormality of the
polynomials $p_{n-1}(x;\mu_{j+1})$ and $p_n(x;\mu_{j})$ and
the reproducing
property of the kernels $K_{n-1}(x,c;\mu_j)$ we have  $$
\left({k_n(\mu_j)
\over k_{n-1}(\mu_{j+1})}\right)^2=1+{p_n(c;\mu_j)^2 \over
K_{n-1}(c,c;\mu_j)}$$ We want to point out that $$\lim_n
{k_n(\mu_j) \over
k_{n-1}(\mu_{j+1})}=1 \qquad {\rm whenever} \qquad \mu_j \in
M(0,1) \quad c
\in [-1,1] \eqno(4)$$ (see [9, Theorem 3 on p.26]), that we
will use later.

Since the polynomials $p_n(x;\mu_{1})$ satisfy
$$K_n(c,c)p_n(x;\mu_{1})={k_n(\mu_1) \over
k_{n+1}}[p_{n+1}'(c)K_n(x,c)-p_{n+1}(c)K_n^{(0,1)}(x,c)]$$ we
can write
$q_n(c)$ and $q_n'(c)$ as follows $$ \eqalign {  q_n(c)
&={\gamma_n \over
k_n D_n} [p_n(c)-N{k_n \over
k_{n-1}(\mu_1)}p_{n-1}'(c;\mu_1)K_{n-1}(c,c)]\cr   q_n'(c)
&={\gamma_n
\over k_n D_n} [p_n'(c)+M{k_n \over
k_{n-1}(\mu_1)}p_{n-1}(c;\mu_1)K_{n-1}(c,c)\}]\cr } \eqno
(5)$$

If we represent the kernels $K_{n-1}(x,c)$ and
$K_{n-1}^{(0,1)}(x,c)$ in
terms of the polynomials $p_n(x)$ and $p_n(x;\mu_j)$ with
$j=1,2$  we can
obtain (see [1, Proposition 2.2])

\proclaim Proposition 1. Let $p_n$ be the orthonormal
polynomials for the
measure $\mu$ and $c \in \R$ such that the condition
$p_n(c)p_{n-1}(c;\mu_1) \not= 0$ is satisfied for  every $n
\in \N$. Then,
the polynomials $q_n$ orthonormal with respect to the inner
product (1)
verify the formula
$$q_n(x)=A_np_n(x)+B_n(x-c)p_{n-1}(x;\mu_1)+C_n(x-c)^2p_{n-2
}(x;\mu_2)
\eqno (6)$$ with $$A_n={\gamma_n\over k_n}(1-\alpha_n) 
\qquad
B_n={\gamma_n \over k_{n-1}(\mu_1)} (\alpha_n-\beta_n) 
\qquad
C_n={\gamma_n \over k_{n-2}(\mu_2)} \beta_n \eqno (6.1)$$
where
$$1-\alpha_n=D_n^{-1} \left[1-N{k_n \over
k_{n-1}(\mu_1)}{p_{n-1}'(c;\mu_1)
\over p_n(c)}K_{n-1}(c,c)\right] \eqno (6.2)$$
$$\beta_n=NK_{n-2}(c,c;\mu_1)D_n^{-1}\left[{k_{n-1}(\mu_1)
\over
k_n}{p_n'(c) \over p_{n-1}(c;\mu_1)}+MK_{n-1}(c,c)\right]
\eqno (6.3)$$ 

\noindent {\bf Remark.} Since all the zeros of the polynomials
$p_n(x)$ and
$p_{n-1}(x;\mu_1)$ are in the interior of the convex hull of
${\rm supp
\enspace \mu}$, then the formula (6) is true  whenever $c$ is
not an
interior point of the convex hull of ${\rm supp \enspace
\mu}$.  \bigskip
From (2), it is obvious that $${\gamma_n \over
k_n}=\int_{\R}q_np_nd\mu
=\langle q_n, p_n \rangle -Mp_n(c)q_n(c)-Np_n'(c)q_n'(c)$$
and then by
straightforward calculations we find, (see [1] or [2]) 
$${\gamma_n \over
k_n}=\left( {D_n \over D_{n+1}} \right)^{1/2} \eqno (7)$$

\bigskip

In the sequel we consider the inner product (1) when the
measure $\mu$ is
the Jacobi weight and $c=1$, that is  $$\langle f,g \rangle
=\int_{[-1,1]}fgw_{\alpha,\beta}dx+M f(1)g(1)+Nf'(1)g'(1)
\eqno (8) $$
where $w_{\alpha,\beta}(x)=(1-x)^{\alpha} (1+x)^{\beta}$
with $\alpha
,\beta >-1$ and $M,N \geq 0$.  
Let $P_n^{(\alpha,\beta)}$ be the Jacobi polynomials with the
normalization
condition \break $P_n^{(\alpha,\beta)}(1)=\displaystyle
{\Gamma(n+\alpha
+1) \over \Gamma(\alpha +1)n!}$ and  $p_n^{(\alpha,\beta)}$ 
the Jacobi
orthonormal polynomials. We denote by $q_n^{(\alpha,\beta)}$
the
polynomials orthonormal with respect to the inner product
$(8)$.

Some basic properties of Jacobi polynomials, (see [11],
Chapter IV), we
will need in the following, are given below. Throughout this
paper we use
the notation $z_n \cong w_n$ when the sequence $z_n/w_n$
converges to 1.

$$P_n^{(\alpha,\beta)}(1) \cong {n^\alpha \over \Gamma(\alpha
+1)} \eqno
(9)$$ $${d \over dx}P_n^{(\alpha,\beta)}(x)={n +\alpha
+\beta +1 \over
2}P_{n-1}^{(\alpha +1, \beta +1)}(x) \eqno (10)$$ $$\Vert
P_n^{(\alpha,\beta)} \Vert^2={2^{\alpha +\beta
+1}\Gamma(n+\alpha
+1)\Gamma(n+\beta +1) \over (2n+\alpha +\beta
+1)n!\Gamma(n+\alpha +\beta
+1)} \cong 2^{\alpha +\beta}n^{-1} \eqno (11)$$
$$a_n={\Gamma(2n+\alpha
+\beta +1) \over 2^n n!\Gamma(n+\alpha +\beta +1)} \cong  
2^{n+\alpha
+\beta}n^{-1/2} \eqno (12)$$ where
$P_n^{(\alpha,\beta)}(x)=a_nx^n+..$.

From (9)-(12), we have for Jacobi orthonormal polynomials:
$$p_n^{(\alpha,\beta)}(1) \cong {n^{\alpha +(1/2)} \over 
2^{(\alpha
+\beta) /2}\Gamma(\alpha +1)}  \eqno (13)$$
$$(p_n^{(\alpha,\beta)})'(1)
\cong {n^{\alpha +(5/2)} \over  2^{(\alpha +\beta +2)/2}
\Gamma(\alpha +2)}
\eqno (14)$$

From these formulas we can deduce

\proclaim Lemma 1. The following estimates hold:
 $$K_n(1,1) \cong {n^{2\alpha +2} \over 2^{\alpha +\beta
+1}\Gamma(\alpha
+1) \Gamma(\alpha +2)} \eqno (15)$$ $$K_n^{(0,1)}(1,1) \cong
{n^{2\alpha
+4} \over 2^{\alpha +\beta +2}\Gamma(\alpha +1) \Gamma(\alpha
+3)} \eqno
(16)$$ $$K_n^{(1,1)}(1,1) \cong {\alpha +2 \over 2^{\alpha
+\beta
+3}\Gamma(\alpha +2) \Gamma(\alpha +4)}n^{2\alpha +6} \eqno
(17)$$

{\bf Proof:} Because of the reproducing property of the
kernels, $K_n(x,1)$
is a polynomial of degree $n$, orthogonal with respect to the
weight
$w_{\alpha +1,\beta}$, that is, for each $n$ there exists a
constant $c_n$
such that  $K_n(x,1)=c_np_n^{(\alpha +1,\beta)}(x)$.
Comparing the leading
coefficients we get $$K_n(x,1)={\Vert P_n^{(\alpha
+1,\beta)} \Vert \over
\Vert P_n^{(\alpha,\beta)} \Vert} {n+\alpha +\beta +1 \over
2n+\alpha
+\beta +1}p_n^{(\alpha,\beta)}(1) p_n^{(\alpha +1,\beta)}(x)
\eqno (18)$$
Now, (15) follows from (11) and (13).

If we derive (18) and evaluating at $x=1$, by using (11),
(13) and (14), we
deduce (16).

To obtain the estimate for $K_n^{(1,1)}(1,1)$ we can consider
the formula
$$K_n(1,1)K_n^{(1,1)}(1,1)-(K_n^{(0,1)}(1,1))^2=
K_{n-1}(1,1;w_{\alpha
+2,\beta})K_n(1,1) \eqno (19)$$ (see [1, Formula $(2.9')$]).
Now (17)
follows from (19), taking into account (15) and (16). \fin

The above lemma and (19) allow us to deduce easily the
asymptotic behaviour
of $D_n$, (see formula (3)).

From now on $C$ will denote a positive constant independent
of $n$, but
possibly different in each ocurrence.

\proclaim Lemma 2. There exists a positive constant $C$ such
that: \item{
a)} if $MN>0$, then $$D_n \cong
MN[K_{n-1}(1,1)K_{n-1}^{(1,1)}(1,1)-(K_{n-1}^{(0,1)}(1,1))^2]
\cong
Cn^{4\alpha +8}$$ \item{ b)} if $M=0$ and $N>0$ then $$D_n
\cong
NK_{n-1}^{(1,1)}(1,1) \cong Cn^{2\alpha +6}$$

Taking in mind (7), a consequence of the above lemma is the
following

\proclaim Corollary 1. Let $k_n$ and $\gamma_n$ be the leading
coefficients
of the polynomials $p_n^{(\alpha,\beta)}$ and
$q_n^{(\alpha,\beta)}$
respectively. Then $\displaystyle \lim_n {\gamma_n \over
k_n}=1$.

\bigskip

{\noindent \bf 3. Estimates for Jacobi-Sobolev polynomials
$q_n^{(\alpha,\beta)}$ on $[-1,1]$ }  \bigskip

In this section, we analyze the behaviour of the Jacobi
Sobolev-type
polynomials $q_n^{(\alpha,\beta)}$ orthonormal with respect to
(8) on
$[-1,1]$.

In order to do this we will estimate the size of  the
coefficients which
appear in their representation in terms of Jacobi polynomials,
see
Proposition 1.

\proclaim Theorem 1. Let $\mu$ be the Jacobi measure, $c=1$
and $A_n$,
$B_n$ and $C_n$ the corresponding coefficients in Proposition
1. There
exists a positive constant $C$ such that:  \item{ a)} if
$MN>0$ then, $A_n
\cong -Cn^{-2\alpha -2} \qquad B_n \cong Cn^{-2\alpha -2}
\qquad C_n \cong
1$  \item{ b)} if $M=0$ and $N>0$ then, $A_n \cong {-1 \over
\alpha +2}
\qquad B_n \cong 1 \qquad C_n \cong {1 \over \alpha +2}$.

{\bf Proof:} Firstly, note that because of (4) and Corollary
1, $\displaystyle {\gamma_n \over k_n}$,  $\displaystyle
{\gamma_n \over
k_{n-1}(w_{\alpha +2,\beta})}$ and  $\displaystyle {\gamma_n
\over
k_{n-2}(w_{\alpha +4,\beta})}$ converge to 1. So, from (6.1),
the
asymptotic behaviour of $A_n$, $B_n$ and $C_n$ only depends on
$\alpha_n$
and $\beta_n$.

a) Assume $MN>0$. Using (13)-(15), we can see that, in formula
(6.2), the
term in  brackets  tends to $-\infty $ like $-n^{2\alpha +6}$.
Since, by
Lemma 2, $D_n \cong Cn^{4\alpha +8}$ it follows that $\alpha_n
\to 1$ and
$A_n \cong -Cn^{-2\alpha -2}$.

Applying formulas (13)-(15) and Lemma 2 in (6.3), we obtain
that $\beta_n
\to 1$; hence $\alpha_n -\beta_n \to 0$. Handling as above, it
is not
difficult to deduce that $B_n \cong Cn^{-2\alpha -2}$.

The result for $C_n$ is immediate.

b) Assume $M=0$ and $N>0$. Lemma 2 and formulas (13)-(15)
lead to  $$
D_n^{-1}{N K_{n-1}(1,1) (p_{n-1}^{(\alpha +2,\beta)})'(1)
\over
p_n^{(\alpha,\beta)}(1)} \to {1 \over \alpha +2}$$  which,
since $D_n \cong
Cn^{2\alpha +6}$, implies that $1-\alpha_n \to -1/(\alpha
+2)$. As to
$\beta_n$, arguing in a similar way we get that $\beta_n \to
1/(\alpha +2)$
and the assertion follows. \fin \bigskip

Now, we can give the asymptotic behaviour of the polynomials
$q_n^{(\alpha,\beta)}$ and $(q_n^{(\alpha,\beta)})'$ at the
ends of the
interval $[-1,1]$ for $M \geq 0$ and $N > 0$.

\bigskip 
\proclaim Theorem 2. There exists a positive constant $C$ such
that the
following estimates

$$q_n^{(\alpha,\beta)}(-1) \cong p_n^{(\alpha,\beta)}(-1)
\cong C(-1)^n
n^{\beta +(1/2)}$$ 
$$(q_n^{(\alpha,\beta)})'(-1) \cong
(p_n^{(\alpha,\beta)})'(-1) \cong
C(-1)^n n^{\beta +(5/2)}$$

$$ q_n^{(\alpha,\beta)}(1) \cong \cases{ -Cn^{-\alpha -(3/2)}
& if $MN>0$ \cr
 -Cn^{\alpha +(1/2)} & if $M=0$, $N>0$ \cr }$$

$$(q_n^{(\alpha,\beta)})'(1) \cong Cn^{-\alpha -(7/2)}$$
\noindent {\it hold.}

{\bf Proof:} Consider the representation of
$q_n^{(\alpha,\beta)}$ in terms
of Jacobi orthonormal polynomials $$q_n^{(\alpha,\beta)}(x)
= A_n
p_n^{(\alpha,\beta)}(x) + B_n(x-1) p_{n-1}^{(\alpha
+2,\beta)}(x)
+C_n(x-1)^2 p_{n-2}^{(\alpha +4,\beta)}(x)$$ Evaluating at
$x=-1$ and
taking into account that $p_n^{(\alpha,\beta)}(-x) = (-1)^n
p_n^{(\beta,\alpha)}(x)$, for all $x \in [-1,1]$, we have

$$q_n^{(\alpha,\beta)}(-1) = (-1)^n [A_n
p_n^{(\alpha,\beta)}(1) + 2B_n
p_{n-1}^{(\beta, \alpha +2)}(1) +4C_n p_{n-2}^{(\beta, \alpha
+4)}(1)]$$
Theorem 1 and (13) yield $\displaystyle \lim_n
{q_n^{(\alpha,\beta)}(-1)
\over p_n^{(\alpha,\beta)}(-1)}=1$, whenever $M \geq 0$ and
$N > 0$.

Deriving in the above expression of $q_n^{(\alpha,\beta)}(x)$
and
proceeding as before, from (13), (14) and Theorem 1,  we
obtain that
$\displaystyle \lim_n {(q_n^{(\alpha,\beta)})'(-1) \over
(p_n^{(\alpha,\beta)})'(-1)}=1$, whenever $M \geq 0$ and $N
> 0$.

To give the asymptotic behaviour at the point $1$, we can use
similar
arguments. However, we want to point out that to estimate
$(q_n^{(\alpha,\beta)})'(1)$ it is easier to apply formula (5)
written for
Jacobi polynomials and $c=1$, that is
$$(q_n^{(\alpha,\beta)})'(1) =
{\gamma_n \over k_n D_n} [(p_n^{(\alpha,\beta)})'(1) + M{k_n
\over
k_{n-1}(w_{\alpha +2,\beta})}p_{n-1}^{(\alpha
+2,\beta)}(1)K_{n-1}(1,1)]$$
  Now it suffices to apply (4), (13)-(15), Lemma 2 and
Corollary 1. \fin

Note that the polynomials orthogonal with respect to the
measure $\mu
+M\delta_1$ are orthogonal with respect to the inner product
$(1)$ with
$c=1$, $M>0$ and $N=0$. Next we summarize for this
situation the main
results of this section: 
\proclaim Lemma 3. Whenever $M>0$ and $N=0$, there exists a
positive
constant $C$ such that, $$D_n \cong MK_{n-1}(1,1) \cong
Cn^{2\alpha +2}$$
$$A_n \cong Cn^{-2\alpha -2} \qquad B_n \cong 1 \qquad
C_n=0$$
$$q_n^{(\alpha,\beta)}(1) \cong Cn^{-\alpha - (3/2)} \qquad
(q_n^{(\alpha,\beta)})'(1) \cong Cn^{\alpha + (5/2)}$$
$$q_n^{(\alpha,\beta)}(-1) \cong p_n^{(\alpha,\beta)}(-1)
\cong C(-1)^n
n^{\beta + (1/2)}$$ $$(q_n^{(\alpha,\beta)})'(-1)
\cong(p_n^{(\alpha,\beta)})'(-1) \cong C(-1)^n n^{\beta +
(5/2)}$$

{\noindent \bf Remark.} Compare the asymptotic behaviour of
$q_n^{(\alpha,\beta)}(1)$ and $q_n^{(\alpha,\beta)}(-1)$ with
the one of
$q_n^{(\alpha,\beta)}(x)$ for $x \in \C \setminus [-1,1]$
which is well
known since Lemma 16 on p.132 in [10] and Theorem 4 in [7]
lead to
$\displaystyle \lim_n { q_n^{(\alpha,\beta)}(x) \over
p_n^{(\alpha,\beta)}(x)} =1$ uniformly for $x$ on compact
sets of  $\C
\setminus [-1,1]$, whenever $M \geq 0$ and $N \geq 0$.
Concerning the
asymptotic behaviour of $p_n^{(\alpha,\beta)}(x)$ out of
$[-1,1]$, see [10,
Theorem 8.21.7].  \bigskip Next we are going to find bounds
for the
polynomials $q_n^{(\alpha,\beta)}$. First, we need to recall a
property
satisfied by Jacobi polynomials.  
Theorem 7.32.2 of [11] shows that there is a constant C
independent of $x$
and $n$ such that $$n^{1/2}|P_n^{(\alpha,\beta)}(x)| \leq
C(1-x+n^{-2})^{-(\alpha /2)-(1/4)} \qquad 0 \leq x \leq 1$$
Using the fact
that $P_n^{(\alpha,\beta)}(x)=(-1)^n P_n^{(\beta,
\alpha)}(-x)$, we have
that the orthonormal Jacobi polynomials satisfy the
estimate
$$|p_n^{(\alpha,\beta)}(x)| \leq C(1-x+n^{-2})^{-(\alpha
/2)-(1/4)}
(1+x+n^{-2})^{-(\beta /2)-(1/4)} \eqno (20)$$ for all $x \in
[-1,1]$ and $n
\geq 1$, with $\alpha ,\beta >-1$. In the sequel $C$ will
denote a positive
constant independent  of $n$ and $x$, but possibly different
in each
ocurrence.

We will find that similar bounds are valid for the polynomials
$q_n^{(\alpha, \beta)}$ with $M,N \geq 0$.

\proclaim Theorem 3. There exists a constant $C$ such that for
each $x \in
[-1,1]$, $n \geq 1$ and $\alpha ,\beta >-1$
$$|q_n^{(\alpha,\beta)}(x)|
\leq C(1-x+n^{-2})^{-(\alpha /2)-(1/4)} (1+x+n^{-2})^{-(\beta
/2)-(1/4)}
\eqno (21)$$

{\bf Proof:} It suffices to prove the result for $n$ large
enough.  
We know, (see Proposition 1), that the polynomials
$q_n^{(\alpha,\beta)}$
satisfy the representation formula

$$q_n^{(\alpha,\beta)}(x)=A_np_n^{(\alpha,\beta)}(x)+B_n(x-1
)p_{n-1}^{(\alpha
+2,\beta)}(x)+C_n(x-1)^2p_{n-2}^{(\alpha +4,\beta)}(x)$$

 Since the coefficients $A_n$, $B_n$ and $C_n$ are bounded
(Theorem 1 and
Lemma 3) and the boundedness (20) for $p_n^{(\alpha,
\beta)}(x)$ is also
true for $(1-x)p_{n-1}^{(\alpha +2,\beta)}(x)$ and $(1-x)^2
p_{n-2}^{(\alpha +4,\beta)}(x)$ for all $x \in [-1,1]$ and $n
\geq 2$, the
statement follows. \fin

As a consequence, whenever $\alpha ,\beta \geq -1/2$, we get a
bound
independent of $n$   $$|q_n^{(\alpha,\beta)}(x)| \leq
C(1-x)^{-(\alpha
/2)-(1/4)} (1+x)^{-(\beta /2)-(1/4)}$$ for all $x \in (-1,1)$.

In particular, if $\alpha=\beta=0$, we have
$|q_n^{(\alpha,\beta)}(x)| \leq
C(1-x^2)^{-(1/4)}$    for all $x \in (-1,1)$. A similar result
has been
obtained in [8] for the polynomials orthonormal with respect
to the inner
product  $\langle f,g \rangle
=\int_{[-1,1]}fgd\mu_0+\int_{[-1,1]}f'g'd\mu_1$ with
$d\mu_0={1 \over 2 }dx
+ M(\delta_1 + \delta_{-1})$ and $d\mu_1=N(\delta_1 +
\delta_{-1})$.

Now, from Theorem 3, we can deduce an upper bound of the
maximum of
$q_n^{(\alpha,\beta)}(x)$ on $[-1,1]$.

\proclaim Corollary 2. There exists a constant $C$ such that
for each $n
\geq 1$ we have

$$ \max_{-1 \leq x \leq 1} |q_n^{(\alpha,\beta)}(x)| \leq
\cases{
Cn^{q+(1/2)} & if $q \geq -1/2$  \cr
 C & if $q \leq -1/2$ \cr } $$
{\noindent \it where $q = \max \{\alpha,\beta \}$.}

{\bf Proof:} The inequalities $1 \leq 1+x+n^{-2} \leq 3$ and
$n^{-2} \leq 1-x+n^{-2} \leq 2$ hold for $x \in [0,1]$. Therefore, from
$(21)$, it
follows that

$$ |q_n^{(\alpha,\beta)}(x)| \leq \cases{ Cn^{\alpha+(1/2)} &
if $\alpha
\geq -1/2$  \cr
 C & if $\alpha \leq -1/2$ \cr } $$

{\noindent for all $x \in [0,1]$.}

A similar argument leads to

$$ |q_n^{(\alpha,\beta)}(x)| \leq \cases{ Cn^{\beta+(1/2)} &
if $\beta \geq
-1/2$  \cr
 C & if $\beta \leq -1/2$ \cr } $$ for all $x \in [-1,0]$. The
assertion
follows easily. \fin

Concerning the asymptotic behaviour of the
$q_n^{(\alpha,\beta)}$ on
$[-1,1]$, by the previous Section we know estimates for these
polynomials
at the end points of the support of the Jacobi weight. What
about the
asymptotic behaviour of the $q_n^{(\alpha,\beta)}$ on $(-1,1)$
?

The Jacobi orthonormal polynomials verify $$p_n^{(\alpha
,\beta)}(x)=r_n^{\alpha ,\beta}(1-x)^{-(\alpha
/2)-(1/4)}(1+x)^{-(\beta
/2)-(1/4)} \cos(k\theta +\gamma)+O(n^{-1}) \eqno (22)$$
$k=n+{\alpha
+\beta+1 \over 2}$, $\gamma = -(\alpha +1) \pi /2$ and 
$r_n^{\alpha
,\beta}=\displaystyle {2^{(\alpha +\beta +1)/2} (\pi
n)^{-1/2} \over \Vert
P_n^{(\alpha ,\beta)} \Vert } \rightarrow \left({2 \over
\pi}\right)^{1/2}$

{\noindent} uniformly for $x$ on compact sets of $(-1,1)$,
(see [11],
Theorem 8.21.8]).

Now, we will show that the polynomials $q_n^{(\alpha ,\beta)}$
have a
similar asymptotic behaviour to the one of $p_n^{(\alpha
,\beta)}$ on the
interval $(-1,1)$.

\proclaim Theorem 4. Let $q_n^{(\alpha,\beta)}$ the
polynomials orthonormal
with respect to (8) and $A_n$, $B_n$ and $C_n$ the
corresponding
coefficients which appear in Proposition 1. Then $$
q_n^{(\alpha,\beta)}(x)=s_n^{\alpha ,\beta}(1-x)^{-(\alpha
/2)-(1/4)}(1+x)^{-(\beta /2)-(1/4)} \cos(k\theta
+\gamma)+O(n^{-1}) $$
$$s_n^{\alpha ,\beta}=A_nr_n^{\alpha
,\beta}+B_nr_{n-1}^{\alpha +2
,\beta}+C_nr_{n-2}^{\alpha +4 ,\beta} \rightarrow \left({2
\over
\pi}\right)^{1/2}$$
 uniformly for $x$ on compact sets of $(-1,1)$. Therefore,
$\displaystyle
\lim_n[ q_n^{(\alpha,\beta)}(x) -
p_n^{(\alpha,\beta)}(x)]=0$ uniformly for $x$ on compact sets of $(-1,1)$.  

{\bf Proof:} By Proposition 1
$$q_n^{(\alpha,\beta)}(x)=A_np_n^{(\alpha,\beta)}(x)+B_n(x-1
)p_{n-1}^{(\alpha
+2,\beta)}(x)+C_n (x-1)^2p_{n-2}^{(\alpha +4,\beta)}(x)$$ From
(22), we have
$$\eqalign{ q_n^{(\alpha,\beta)}(x) &=(1-x)^{-(\alpha
/2)-(1/4)}(1+x)^{-(\beta /2)-(1/4)} \cos(k\theta
+\gamma)[A_nr_n^{\alpha
,\beta}+B_nr_{n-1}^{\alpha +2 ,\beta}+C_nr_{n-2}^{\alpha +4
,\beta}] \cr
&+[A_n+B_n(x-1)+C_n(x-1)^2]+O(n^{-1}) \cr }$$ {\noindent}
uniformly for $x$
on compact sets of $(-1,1)$.

Since the asymptotic behaviour of the coefficients $A_n$,
$B_n$ and $C_n$
obtained in the previous section, we get $$
q_n^{(\alpha,\beta)}(x)=s_n^{\alpha ,\beta}(1-x)^{-(\alpha
/2)-(1/4)}(1+x)^{-(\beta /2)-(1/4)} \cos(k\theta
+\gamma)+O(n^{-1}) $$
{\noindent}and $\displaystyle \lim_n s_n^{\alpha
,\beta}=\left({2 \over
\pi}\right)^{1/2}$.

Therefore $\displaystyle
q_n^{(\alpha,\beta)}(x)={s_n^{\alpha ,\beta} \over
r_n^{\alpha ,\beta}}p_n^{(\alpha,\beta)}(x)+O(n^{-1})$ and we
can write
$$q_n^{(\alpha,\beta)}(x)
-p_n^{(\alpha,\beta)}(x)=\left({s_n^{\alpha
,\beta} \over r_n^{\alpha
,\beta}}-1\right)p_n^{(\alpha,\beta)}(x)+O(n^{-1})$$
{\noindent} uniformly
for $x$ on compact sets of $(-1,1)$. Thus the result follows.
\fin

{\noindent \bf Remark.} From (20) we have $\vert
p_n^{(\alpha,\beta)}(x)
\vert \leq C$ for $x$ on compact sets of $(-1,1)$. Then
$\displaystyle
\lim_n[q_n^{(\alpha,\beta)}-p_n^{(\alpha,\beta)}]=0$
uniformly on compact
sets of $(-1,1)$ could be also deduced applying Theorem 5 in
[7] and
formula (10) of Lemma 16 in [10].
\bigskip {\noindent \bf 4. Estimates for the kernels} \bigskip

It is known, (Nevai [10, Lemma 5 on p. 108]), that the kernels
associated
with Jacobi polynomials satisfy the estimate  $$K_n(x,x) \sim
n(1-x+n^{-2})^{-\alpha -(1/2)}(1+x+n^{-2})^{-\beta -(1/2)}
\eqno (23)$$
\noindent uniformly in $\vert x \vert \leq 1$, $n \geq 1$,
where by $f_n(x)
\sim g_n(x)$ we mean that there exist some positive constants
$C_1$ and
$C_2$ such that $C_1f_n(x) \leq g_n(x) \leq C_2f_n(x)$ for all
$x \in
[-1,1]$ and $n \in \N$.

We want to find similar estimates for the new kernels.

Let $L_n(x,y)$ be the kernels relative to the inner product
(8). If we
consider their expansion in terms of Jacobi orthonormal
polynomials, we can
deduce, (see [1, p.744]),  $$L_n(x,y) = K_n(x,y) -M
L_n(y,1)K_n(x,1) - N
L_n^{(0,1)}(y,1) K_n^{(0,1)}(x,1) \eqno (24)$$  \vfill  \eject
with
$$L_n(x,1)= D_{n+1}^{-1} \left([1+NK_n^{(1,1)}(1,1)]K_n(x,1)
-
NK_n^{(0,1)}(1,1)K_n^{(0,1)}(x,1) \right)$$

$$L_n^{(0,1)}(x,1)= D_{n+1}^{-1}
\left([1+MK_n(1,1)]K_n^{(0,1)}(x,1) - M
K_n^{(0,1)}(1,1)K_n(x,1) \right)$$ Inserting $L_n(x,1)$ and
$L_n^{(0,1)}(x,1)$ in (24) and taking $y=x$, we get
$$\eqalign{ L_n(x,x) &
= K_n(x,x) - D_{n+1}^{-1} [ M \{ 1+NK_n^{(1,1)}(1,1) \}
K_n(x,1)^2 \cr  & -
2MN K_n^{(0,1)}(1,1) K_n(x,1) K_n^{(0,1)}(x,1) + N \{
1+MK_n(1,1) \}
K_n^{(0,1)}(x,1)^2] \cr} \eqno (25)$$ If, as usual, we define
the
Christoffel function  $$\Lambda_n(x)=\min \{ \langle p, p
\rangle; \deg p
\leq n, p(x)=1 \}$$  {\noindent}it is easy to see that
$\Lambda_n(x)=[L_n(x,x)]^{-1}$.

We will use the representation (25) to obtain some bounds for
$L_n(x,x)$.

\proclaim Theorem 5. Let $(L_n(x,y))$ be the kernels relative
to the
polynomials $q_n^{(\alpha, \beta)}$. Then  there exists a
constant $C$ such
that for each $x \in [-1,1]$ and $n \geq1$

$$|L_n(x,x)| \leq Cn(1-x+n^{-2})^{-\alpha -(1/2)}
(1+x+n^{-2})^{-\beta
-(1/2)}$$
{\bf Proof:} From $(23)$ we have for each $x \in [-1,1]$, $n
\geq 1$ and
$\alpha ,\beta >-1$
$$|K_n(x,x)| \leq Cn(1-x+n^{-2})^{-\alpha -(1/2)}
(1+x+n^{-2})^{-\beta
-(1/2)} \eqno (26)$$ Moreover, from (11), (18) and (20), 
$$\eqalign{
|K_n(x,1)| & \leq C |p_n^{(\alpha,\beta)}(1)| |p_n^{(\alpha
+1,\beta)}(x)|
\cr & \leq C n^{\alpha +(1/2)} (1-x+n^{-2})^{-(\alpha
/2)-(3/4)}
(1+x+n^{-2})^{-(\beta /2)-(1/4)} \cr} \eqno (27)$$ for all $x
\in [-1,1]$.

To find a bound for $K_n^{(0,1)}(x,1)$, we will use the
formula
$$K_n^{(0,1)}(x,1) = (x-1) K_{n-1}(x,1;w_{\alpha +2,\beta})
+
{K_n^{(0,1)}(1,1) \over K_n(1,1)}  K_n(x,1) \eqno (28) $$ (see
[1, Formula
(2.9)]), from which, using (27) and Lemma 1, it follows that
$$|K_n^{(0,1)}(x,1)| \leq C n^{\alpha +(5/2)}
(1-x+n^{-2})^{-(\alpha
/2)-(3/4)} (1+x+n^{-2})^{-(\beta /2)-(1/4)} \eqno (29)$$ for
all $x \in
[-1,1]$.

Now it suffices to remind that by Lemmas 1 and 2, whenever
$MN>0$ $$M
D_{n+1}^{-1} [1+NK_n^{(1,1)}(1,1)] \leq C n^{-2\alpha -2}$$
$$2MN
D_{n+1}^{-1} K_n^{(0,1)}(1,1) \leq C n^{-2\alpha -4}$$ $$N
D_{n+1}^{-1}
[1+MK_n(1,1)] \leq C n^{-2\alpha -6}$$  and to observe that
for each $x \in
[-1,1]$, the inequality $n^{-1}(1-x+n^{-2})^{-1} \leq Cn$
holds. For the
other values of the parameters $M$ and $N$, we proceed in a
similar way.
Thus, the result follows. \fin

\bigskip This result gives us only upper bounds. Now we want
to estimate
more accurately $L_n(x,x)$. First, we observe the behaviour of
$L_n(x,x)$
at the end points of the interval $[-1,1]$. Evaluating at
$x=1$ the
expression of $L_n(x,1)$ given in (24) and using (19), we get
$$L_n(1,1)=D_{n+1}^{-1}[1+NK_{n-1}(1,1;w_{\alpha
+2,\beta})]K_n(1,1)$$
Then, the kernels $L_n(1,1)$ are bounded if $M>0$, $N \geq 0$
while
$L_n(1,1) \cong CK_n(1,1)$ if $M=0, N\geq 0$. Note that
the boundedness of
$L_n(1,1)$ depends on the addition of a mass at 1 and not of
the term
involving derivatives.

Moreover from the expression of $L_n(1,1)$ we can recover the
mass $M$.
Indeed, by using Lemmas 1, 2 and 3 it follows that, when $M>0$
and $N \geq
0$, $\displaystyle \lim_n \Lambda_n(1)=M$. Otherwise,  the
mass $N$ can be
recovered from $L_n^{(1,1)}(1,1)$ since, when $M \geq 0$ and
$N >0$,
$\displaystyle \lim_n[L_n^{(1,1)}(1,1)]^{-1}=N$.

As to $L_n(-1,-1)$, it suffices to take $x=y=-1$ in (24)
and we obtain
$$L_n(-1,-1) \cong CK_n(-1,-1) \cong Cn^{2 \beta +2}$$ Next,
we are going
to find uniform estimates for the kernels. When $M>0$,
$L_n(1,1)$ is
bounded, so we give uniform estimates on compact sets not
containing the
mass point $1$.

\proclaim Theorem 6.  a) Suppose $M>0, N \geq 0$. Let
$\varepsilon >0$, then
$$L_n(x,x) \sim n(1-x+n^{-2})^{-\alpha
-(1/2)}(1+x+n^{-2})^{-\beta -(1/2)}
$$ uniformly on $[-1,1-\varepsilon]$,  $n \geq 1$.

{\noindent}{\sl  b) Suppose $M=0, N \geq 0$. Then $$L_n(x,x)
\sim
n(1-x+n^{-2})^{-\alpha -(1/2)}(1+x+n^{-2})^{-\beta -(1/2)} $$
uniformly on
$\vert x \vert \leq 1$, $n \in \N$.}

{\bf Proof:} Because of Theorem 5, it suffices to prove that,
for $n$ large
enough $$L_n(x,x) \geq  Cn(1-x+n^{-2})^{-\alpha
-(1/2)}(1+x+n^{-2})^{-\beta
-(1/2)} $$ uniformly on $[-1,1-\varepsilon]$ when $M >0$ and
on $[-1,1]$
when $M = 0$.

For the sake of simplicity, we write 
$$d(x,n)=n(1-x+n^{-2})^{-\alpha
-(1/2)}(1+x+n^{-2})^{-\beta -(1/2)}$$

a) Let $N>0$. Using Lemmas 1 and 2 and formulas (27) and (29),
we obtain
that the three last summands in (25) are bounded by
$Cd(x,n)n^{-2}(1-x+n^{-2})^{-1}$. Thus, taking into account
(23), the
result follows. For $N=0$, we handle in a similar way.

b) For $N=0$ the result is obvious because of
$L_n(x,x)=K_n(x,x)$. Suppose
$N>0$, as $D_{n+1}=1+NK_n^{(1,1)}(1,1)$, from (25) we have
$$L_n(x,x) \geq
ND_{n+1}^{-1}[K_n^{(1,1)}(1,1)K_n(x,x)-K_n^{(0,1)}(x,1)^2]$$
and using,
again, the estimates for the kernels and (28) we can deduce
the result. \fin

Now we consider the analogue of the Szeg\"o extremum problem
for the inner
product (8).

It is known that the generalized Szeg\"o extremum problem,
associated with
a finite positive Borel measure on the real line, consists of
finding
$\displaystyle \lim_n \lambda_n(x;\mu)$ with
$\lambda_n(x;\mu)$ the
Christoffel functions corresponding to $\mu$. A solution of
this problem,
when $\mu$ belongs to the Szeg\"o class of the interval
$[-1,1]$, has been
given in [9, Theorem 5] by proving that $\displaystyle \lim_n
n\lambda_n(x;\mu)=\pi \mu'(x)(1-x^2)^{1/2}$ for almost every
$x \in
[-1,1]$, where $\mu'$ is almost everywhere the Radon-Nikodym
derivative of
$\mu$, (see [9]).

\proclaim Theorem 7. Let $\Lambda_n$ the Christoffel functions
associated
with (8). Then $$ \lim_n n\Lambda_n(x)=\pi w_{\alpha
,\beta}(x)(1-x^2)^{1/2}$$ for almost every $x \in [-1,1]$.

{\bf Proof:} Because of M\'at\'e-Nevai-Totik result, above
quoted, we only
need to prove $\displaystyle \lim_n n^{-1}L_n(x,x)=\lim_n
n^{-1}K_n(x,x)$,
$x \in [-1,1]$. Thus, by (25), it suffices to deduce 
$$\eqalign{ \lim_n &
D_{n+1}^{-1}[1+NK_n^{(1,1)}(1,1)]K_n(x,1)^2=0 \cr \lim_n &
D_{n+1}^{-1}K_n^{(0,1)}(1,1)K_n(x,1)K_n^{(0,1)}(x,1)=0 \cr
\lim_n &
D_{n+1}^{-1}[1+MK_n(1,1)]K_n^{(0,1)}(x,1)^2=0 \cr }$$
{\noindent} for every
$x \in (-1,1)$ and this follows by considering (18), (13),
(22) and (28).
\fin  \bigskip From the results of Section 2 and formulas (27)
and (29),
the following bounds for $L_n(x,1)$ and $L_n^{(0,1)}(x,1)$ can
also be
obtained: \medskip

{\noindent}{\bf Theorem 8.} {\sl There exists a constant $C$
such that for
each $x \in [-1,1]$ and $n \geq 1$} $$\matrix{ &\vert L_n(x,1)
\vert \leq
C(1+x+n^{-2})^{-(\beta/2)-(1/4)} \qquad &{\rm if}\qquad M>0
\cr &\vert
L_n(x,1) \vert \leq Cn^{2\alpha
+4}(1+x+n^{-2})^{-(\beta/2)-(1/4)} \qquad
&{\rm if}\qquad M=0 \cr &\vert L_n^{(0,1)}(x,1) \vert  \leq
C(1+x+n^{-2})^{-(\beta/2)-(1/4)} \qquad &{\rm if}\qquad  N>0
\cr &\vert
L_n^{(0,1)}(x,1) \vert  \leq Cn^{2\alpha
+4}(1+x+n^{-2})^{-(\beta/2)-(1/4)}
\qquad &{\rm if}\qquad N=0 \cr}$$

Notice that the above bounds for $L_n(x,1)$ when $M>0$ and for
$L_n^{(0,1)}(x,1)$ when $N>0$ are, respectively, sharper than
the ones for $K_n(x,1)$ and $K_n^{(0,1)}(x,1)$ (see formulas (27) and
(29)).

\medskip  {\noindent \bf Remark}. Some of the above results about the
kernels  appear
in [5] for $w$ a generalized Jacobi weight and $N=0$. 
\bigskip

Finally, it
is worth observing that if in the product (1) $\mu$ is the
Jacobi measure
and we take $c=-1$, since Jacobi polynomials satisfy
$p_n^{(\alpha,\beta)}(-1) = (-1)^n p_n^{(\beta,\alpha)}(1)$,
we get the
same results as above but exchanging $\alpha$ and $\beta$.

\bigskip {\noindent \bf References}

\medskip [1] M. Alfaro, F. Marcell\'an, M.L. Rezola and A.
Ronveaux, On
orthogonal polynomials of Sobolev type: Algebraic properties
and zeros, {\it SIAM J. Math. Anal. 23 (1992) 737-757.}

[2] M. Alfaro, F. Marcell\'an, M.L. Rezola and A. Ronveaux,
Sobolev-type
orthogonal polynomials: The nondiagonal case, {\it J. Approx.
Theory 83
(1995) 266-287.} 

[3] H. Bavinck and H.G. Meijer, Orthogonal polynomials with
respect to a
symmetric inner product involving derivatives, {\it Applicable
Analysis 33
(1989) 103-117.}

[4] H. Bavinck and H.G. Meijer, On orthogonal polynomials with
respect to
an inner product involving derivatives: zeros and recurrence
relations,
{\it Indag. Math. (N.S.)
 1 (1990) 7-14.}

[5] J.J. Guadalupe, M. P\'erez, F.J. Ruiz and J.L. Varona,
Asymptotic
behaviour of orthogonal polynomials relative to measures with
mass points,
{\it Mathematika 40 (1993) 331-344.}

[6] G. L\'opez, F. Marcell\'an and W. Van Assche, Relative 
asymptotics for
polynomials orthogonal with respect to a discrete Sobolev 
inner product,
{\it Constr. Approx. 11 (1995) 107-137.}

[7] F. Marcell\'an and W. Van Assche, Relative asymptotics for
orthogonal
polynomials with a Sobolev inner product, {\it J. Approx.
Theory 72 (1992)
192-209.}

[8] F. Marcell\'an and B. Osilenker,  Estimates for
polynomials orthogonal
with respect to some Legendre-Sobolev type inner product.
Submitted

[9] A. M\'at\'e, P.G. Nevai and V. Totik, Szeg\"o's extremum
problem on the
unit circle, {\it Ann. Math. 134 (1991) 433-453.}

[10] P.G. Nevai, {\it Orthogonal Polynomials}, Memoirs Amer.
Math. Soc.
213,  Amer. Math. Soc., Providence, RI, 1979.

[11] G. Szeg\"o, {\it Orthogonal Polynomials}, Amer. Math.
Soc. Colloq.
Pub. 23,  Amer. Math. Soc., Providence, RI, 1975 (4th
edition).

\bye